\documentclass[a4paper,10pt,conference]{IEEEconf}
\usepackage{cite}
\usepackage{amsmath}
\usepackage{amssymb}
\usepackage{graphicx}
\usepackage{psfrag}
\begin{document}
\newtheorem{theorem}{Theorem}[section]
\newtheorem{definition}{Definition}[section]

\title {On the Lossless Property of a Class of Nonlinear Quantum Systems}
\author{Aline I. Maalouf and Ian R. Petersen \\
School of Engineering and Information Technology,\\ University of New South Wales
at the Australian Defence Force
Academy,\\ Canberra, ACT 2600 \tt\small a.maalouf@adfa.edu.au \tt\small i.r.petersen@gmail.com
\thanks{This work was completed with the support of a University
of New South Wales Postgraduate Award and the Australian Research Council.}}
\maketitle
\thispagestyle{empty}
\pagestyle{empty}
\maketitle
\begin{abstract}
In this paper, the physical realizability condition for a specific class of nonlinear quantum systems is related to the lossless property of nonlinear dissipative systems having a specific storage function. %The physical realizability property for this class of nonlinear quantum systems determines whether a given set of nonlinear quantum stochastic differential equations corresponds to a physical nonlinear quantum system satisfying the laws of quantum mechanics.
\end{abstract}
\section{Introduction}
%In the physics literature, methods have been developed to model a wide range of open quantum systems, such as those encountered in quantum optics, within the framework of quantum stochastic differential equations \cite{Gardiner85}, \cite{Gardiner00}, \cite{Hudson84}. In this context, quantum noise is used to represent the influence of large heat baths and boson fields, including optical and phonon fields, from which completely positive maps, Lindblad generators and master equations are obtained by taking expectations \cite{James10}.
  There are two types of quantum stochastic differential equations, namely, linear and nonlinear. Linear quantum systems can be described by linear quantum stochastic differential equations and arise mostly in the area of quantum optics; e.g., see \cite{Walls08}, \cite{Gardiner00} and \cite{Bachor04}.
An important class of linear quantum stochastic models describe the Heisenberg evolution of the annihilation and creation operators of several independent open quantum harmonic oscillators that are coupled to external coherent bosonic fields, such as coherent laser beams; e.g., see \cite{Wiseman10}, \cite{Walls08} and \cite{Gardiner00}. A special class of these linear quantum systems is driven by quantum Wiener processes as in \cite{James08} where physical realizability conditions are developed to determine when the linear quantum system under consideration can be regarded as a representation of a linear quantum harmonic oscillator. This imposes some restrictions on the matrices describing the linear quantum model. This notion of physical realizability has been further investigated in \cite{Shaiju09} where the authors prove the equivalence between the algebraic conditions for physical realizability obtained in \cite{James08} and a frequency domain condition that an associated linear system is $(J, J)$-unitary. In fact, \cite{Shaiju09} extends the frequency domain physical realizability results of \cite{Maalouf09}, \cite{Maalouf11-1} and \cite{Maalouf11-2} which apply to linear quantum systems described purely in terms of annihilation operators. More explicitly, in \cite{Maalouf09}, \cite{Maalouf11-1} and \cite{Maalouf11-2}, the physical realizability conditions developed for annihilation-operator linear quantum systems have been related to the lossless property of linear systems and \cite{Shaiju09} generalizes this result to relate the corresponding physical realizability conditions of annihilation and creation operators linear quantum systems to the $(J, J)$-unitary property.

%The question of addressing the physical realizability conditions for nonlinear quantum systems described by nonlinear quantum stochastic differential equations is still open. 

In this paper, we restrict attention to a class of nonlinear quantum systems which is a generalization of the annihilation operator only linear quantum systems considered in \cite{Maalouf09}, \cite{Maalouf11-1} and \cite{Maalouf11-2}. Then, we relate the resulting physical realizability conditions to the lossless property of nonlinear dissipative systems having a specific type of storage function.

%The paper is organized as follows: Section II defines the class of nonlinear quantum systems under consideration. Section III develops conditions under which the commutation relations for the nonlinear quantum systems are preserved. Section IV defines physical realizability in terms of nonlinear open quantum harmonic oscillators. Section V defines the lossless properties of the class of nonlinear quantum systems under consideration. Section VI relates the physical realizability conditions obtained to the lossless property of nonlinear dissipative quantum systems having a specific storage function. Section VII provides an example to illustrate the theory developed. Section VIII concludes the paper.

\section{The Class of Nonlinear Quantum Systems}

We consider an open quantum system G with physical variable space $\mathcal{A}_G$ consisting of operators defined on an underlying Hilbert space $\mathcal{H}_G$. The self-energy of this system is described by a Hamiltonian $\mathcal{H}\in \mathcal{A}_G$. The system is driven by $m$ field channels given by the quantum stochastic process $W$

\begin{equation}\nonumber
W= \left( {\begin{array}{*{20}{c}}
{\begin{array}{*{20}{c}}
{{W_1}}\\
 \vdots \\
{{W_m}}
\end{array}}
\end{array}} \right).
\end{equation}
These describe the annihilation of photons in the field channels and are operators on a Hilbert space $F$, with associated variable space $\mathcal{F}$. In that case, $F$ is the Hilbert space defining an indefinite number of quanta (called a Fock space \cite{Parthasarathy92}), and $\mathcal{F}$ is the space of operators over this space.

We assume the process $W$ is canonical, meaning that we have the following second order Ito products:
\begin{eqnarray}\nonumber
dW_k(t)dW_l(t)^*&=&\delta_{kl}dt;\\\nonumber
dW_k(t)^*dW_l(t)&=&0;\\\nonumber
dW_k(t)dW_l(t)&=&0;\\\nonumber
dW_k(t)^*dW_l(t)^*&=&0
\end{eqnarray}
where $W_k(t)^*$ is the operator adjoint of $W_k(t)$ which is defined on the same Fock space.

The system is coupled to the field through a scattering matrix $S=I$ and a coupling vector of operators $L$ given by

\begin{equation}\nonumber
L=\left( {\begin{array}{*{20}{c}}
{\begin{array}{*{20}{c}}
{{L_1}}\\
 \vdots \\
{{L_m}}
\end{array}}
\end{array}} \right)
\end{equation}
where $L_j \in \mathcal{A}_G$.

The notation $G=(S,L,H)$ is used to indicate an open quantum system specified by the parameters $S, L$ and $H$ where $H$ represents the Hamiltonian of the system. For $S=I$, the Schrodinger equation is given by
\begin{equation}\label{schrodinger}
dU(t)= \left\{dW^\dagger L - L^\dagger dW -\frac{1}{2}L^\dagger L dt-iHdt\right\}U(t)
\end{equation}
with the initial condition $U(0)=I$. Equation (\ref{schrodinger}) determines the unitary motion of the system in accordance of the fundamental laws of quantum mechanics. Note that the notation $^\dagger$ refers to the Hilbert space adjoint.

Given a system annihilation operator $a_l \in \mathcal{A}_G$, its Heisenberg evolution is defined by $a_l(t)=j_t(a_l)=U(t)^* a_l U(t)$ and satisfies
\begin{eqnarray}\label{annihi_eq}
da_l(t)&=&\left(\mathcal{L}_L(t)(a_l(t))-i\left[a_l(t),H(t)\right]\right)\\\nonumber
&& + \left[L(t)^\dagger,a_l(t)\right]dW(t)
\end{eqnarray}
where the notation $\left[A,B\right]=AB-BA$ denotes the commutator of two operators $A$ and $B$.

In equation (\ref{annihi_eq}), all operators evolve unitarily and the notation $\mathcal{L}_L(a_l)$ refers to
\begin{equation}\label{gen_l}
\mathcal{L}_L(a_l)= \frac{1}{2}L^\dagger \left[a_l,L\right]+ \frac{1}{2} \left[L^\dagger,a_l\right]L.
\end{equation}
In this paper, we restrict attention to annihilation only coupling operators; i.e, the coupling operator $L(t)$ depends only on the annihilation operator $a_l(t)$ and not the creation operator $a_l(t)^*$. Therefore, $L(t)$ satisfies the following commutation property: $\left[a_l,L\right]=0$ and then
\begin{equation}\label{gen_l_a}
\mathcal{L}_L(a_l)= \frac{1}{2} \left[L^\dagger,a_l\right]L.
\end{equation}
The generator of the system $G$ is given by
\begin{equation}\label{gen_g}
\mathcal{G}_G(a_l)=-i\left[a_l,H\right]+\mathcal{L}_L(a_l).
\end{equation}
For the case of having $n$ annihilation operators $a_1,\ldots,a_n$, we define
\begin{equation}\nonumber
a= \left( {\begin{array}{*{20}{c}}
{\begin{array}{*{20}{c}}
{{a_1}}\\
 \vdots \\
{{a_n}}
\end{array}}
\end{array}} \right),
\end{equation}
\begin{equation}\label{genl_i}
\mathcal{L}_{{L}_i}(a)=  \frac{1}{2} \left[L^\dagger,a_i\right]L
\end{equation}
and
\begin{equation}\label{geng_i}
\mathcal{G}_{G_i}(a)=-i\left[a_i,H\right]+\mathcal{L}_{L_i}(a).
\end{equation}
Therefore, we can write
\begin{eqnarray}\nonumber \label{annihieq_i}
da_i(t)&=&\left(\mathcal{L}_{L_i}(t)(a_i(t))-i\left[a_i(t),H(t)\right]\right)\\\nonumber
&& + \left[L(t)^\dagger,a_i(t)\right]dW(t)
\end{eqnarray}
and
\begin{eqnarray}\nonumber \label{annihieq_a}
da(t)&=&\left(\mathcal{L}_{L}(t)(a(t))-i\left[a(t),H(t)\right]\right)\\\nonumber
&& + \left[L(t)^\dagger,a(t)\right]dW(t)
\end{eqnarray}
where
\begin{equation}\label{L_a}
\mathcal{L}_{L}(a)= \left( {\begin{array}{*{20}{c}}
{\begin{array}{*{20}{c}}
{{\mathcal{L}_{{L}_1}(a)}}\\
 \vdots \\
{{\mathcal{L}_{{L}_n}(a)}}
\end{array}}
\end{array}} \right)
\end{equation}
and
\begin{equation}\label{G_a}
\left[a,H\right]= \left( {\begin{array}{*{20}{c}}
{\begin{array}{*{20}{c}}
{{\left[a_1,H\right]}}\\
 \vdots \\
{{\left[a_n,H\right]}}
\end{array}}
\end{array}} \right).
\end{equation}
Let
\begin{eqnarray}\nonumber\label{At}
A(a(t), a(t)^\dagger)&=&\left({\begin{array}{*{20}{c}}
{\begin{array}{*{20}{c}}
{{A_1(a(t),a(t)^\dagger )}}\\
 \vdots \\
{{A_n(a(t), a(t)^\dagger)}}
\end{array}}
\end{array}}\right)\\\nonumber&=& \left( {\begin{array}{*{20}{c}}
{\begin{array}{*{20}{c}}
{{\mathcal{L}_{{L}_1}(a)}}\\
 \vdots \\
{{\mathcal{L}_{{L}_n}(a)}}
\end{array}}
\end{array}} \right)-i\left( {\begin{array}{*{20}{c}}
{\begin{array}{*{20}{c}}
{{\left[a_1,H\right]}}\\
 \vdots \\
{{\left[a_n,H\right]}}
\end{array}}
\end{array}} \right)\\\nonumber &=& \mathcal{L}_{L}(a)-i\left[a,H\right]\\ &=&\mathcal{G}_G(a)
\end{eqnarray}
and
\begin{eqnarray}\label{Bt}
B(a(t),a(t)^\dagger)&=&\left({\begin{array}{*{20}{c}}
{\begin{array}{*{20}{c}}
{{B_1(a(t),a(t)^\dagger)}}\\
 \vdots \\
{{B_n(a(t),a(t)^\dagger)}}
\end{array}}
\end{array}}\right)\\\nonumber &=& \left({\begin{array}{*{20}{c}}
{\begin{array}{*{20}{c}}
{{\left[L^\dagger, a_1\right]}}\\
 \vdots \\
{{\left[L^\dagger, a_n\right]}}
\end{array}}
\end{array}}\right)\\&=& \left[L^\dagger, a\right].
\end{eqnarray}
Hence,
\begin{eqnarray}\nonumber
da(t)&=&A(a(t),a(t)^\dagger)dt+B(a(t),a(t)^\dagger)dW(t)\\\nonumber
da(t)^*&=&A(a(t),a(t)^\dagger)^*dt+B(a(t),a(t)^\dagger)^*dW(t)^*.
\end{eqnarray}
Note that for the case of matrices, the notations $^*$ and $^\dagger$ refer respectively to the complex conjugate and the complex conjugate transpose of the matrix in question.
\begin{equation}\label{at}\end{equation}
On the other hand, the components of the output fields are defined by $y(t)=j_t(W(t))=U(t)^*W(t)U(t)$ and satisfy the nonlinear quantum stochastic differential equations
\begin{eqnarray}\nonumber
dy(t)&=&C(a(t))dt+D(t)dW(t)\\\nonumber
dy(t)^*&=&C(a(t))^*dt+D(t)^*dW(t)^*
\end{eqnarray}
\begin{equation}\label{yt}\end{equation}
where $C(a(t))=L(t)$, $C(a(t))^*=L(t)^*$, $D(t)=I$ and $D(t)^*=I$.
Hence, the system $G$ can be described by the following nonlinear quantum stochastic differential equations
\begin{eqnarray}\nonumber\label{sys_a}
d\bar a(t)&=&\bar A(a(t),a(t)^\dagger)dt+\bar B(a(t),a(t)^\dagger)d\bar W(t);\\
d\bar y(t)&=&\bar C(a(t),a(t)^\dagger)dt+\bar D(a(t),a(t)^\dagger)d\bar W(t)
\end{eqnarray}
where $\bar a (t)=\left[ {\begin{array}{*{20}{c}}
{a}\\
{a^*}
\end{array}} \right]$, $\bar A(a,{a^\dag }) = \left[ {\begin{array}{*{20}{c}}
{A(a,{a^\dag })}\\
{A{{(a,{a^\dag })}^*}}
\end{array}} \right]$, $\bar C(a,{a^\dag }) = \left[ {\begin{array}{*{20}{c}}
{C(a)}\\
{C{{(a)}^*}}
\end{array}} \right]$, $\bar B(a,{a^\dag }) = \left[ {\begin{array}{*{20}{c}}
{B(a,a^\dagger)}&0\\
0&{B{{(a,a^\dagger)}^*}}
\end{array}} \right]$, $\bar D(t) = \left[ {\begin{array}{*{20}{c}}
{D(t)}&0\\
0&{D{{(t)}^*}}
\end{array}} \right]$ and $d\bar W(t)=\left[ {\begin{array}{*{20}{c}}
{dW(t)}\\
{dW(t)^*}
\end{array}} \right]$.
\begin{definition}\label{def4}
The class of nonlinear quantum system we consider in this paper are nonlinear quantum systems that can be represented by the QSDES (\ref{sys_a}) such that the matrices $A(a,a^\dagger), B(a,a^\dagger)$ and $ C(a)$ satisfy $\left[B(a,a^\dagger), a^T\right]=0$, $\left[C(a), a^T\right]=0$, $\left[A(a,a^\dagger), a^T\right]=-\left[a, A(a,a^\dagger)^T\right]$ and
\begin{eqnarray}\nonumber \label{prop_matrices}
A_i(a,a^\dagger)&=& \sum\limits_{{k_i} = 0}^{{m_{{k_i}}}} {\sum\limits_{{h_i} = 0}^{{m_{{h_i}}}} {\sum\limits_{l = 1}^n {\sum\limits_{p = 1}^n {{\alpha _{{plk_ih_i}} } }}}} a_p^{{k_i}}{(a_l^*)^{{h_i}}}\\
C_v(a)&=& \sum\limits_{{k_v} = 0}^{{m_{{k_v}}}}\sum\limits_{p = 1}^n \beta_{k_vp}a_p^{k_v}
\end{eqnarray}
for $i=1,\ldots,n$ and $v=1,\ldots, m$ with $m_{k_i}$, $m_{h_i}$, $m_{k_v}$, $\beta_{k_vp}$ and $\alpha _{{plk_i}{h_i}}$ integers.

Moreover, the matrices $A(a,a^\dagger)$, $B(a,a^\dagger)$ and $C(a)$ satisfy the following property:
\begin{eqnarray}\nonumber \label{prop_matrices_2}
&&\frac{1}{{{2{\bar n}}}}\left[ {\bar A{{(a,{a^\dag })}^\dag }{{\bar \theta }^{ - 1}}\bar a,\bar a} \right] - \frac{1}{{{2{\bar n}}}}\left[ {{{\bar a}^\dag }{{\bar \theta }^{ - 1}}\bar A(a,{a^\dag }),\bar a} \right]=\\&& \qquad \bar A(a,{a^\dag }) - \frac{1}{2}\bar B(a,{a^\dag })\bar C(a,a^\dagger)
\end{eqnarray}
where $\bar n= {\left. {\mathop {\sup }\limits_{1 \le i \le n} \left( {{k_i} + {h_i}} \right)} \right|_{{\alpha _{{plk_ih_i}}}\ne 0}}$, $\bar \Theta= \left[ {\begin{array}{*{20}{c}}
\Theta &0\\
0&{{\Theta ^*}}
\end{array}} \right]$ and $\Theta = \left[a,a^\dagger\right]$. For more information, on $\Theta$, refer to Section III below.
\end{definition}

In the sequel, we assume that $dW(t)$ admits the following
decomposition:
\begin{equation}\nonumber
 dW (t) = \beta _{w} (t)dt + d\tilde w(t)
\end{equation}
where $\tilde w(t)$ is the noise part of $W(t)$ and $\beta_{w}(t)$
is an adapted process (see \cite{Parthasarathy92} and \cite{Hudson84}).
The noise $\tilde w(t)$ is a vector of quantum Weiner processes with Ito table
\begin{equation}\label{wIto}
d\tilde w(t)d\tilde w(t)^\dagger = F_{\tilde{w}} dt
\end{equation}
(see \cite{Parthasarathy92}) where $F_{\tilde{w}}$ is a non-negative definite Hermitian matrix. Here, the notation $^\dagger$ represents the adjoint transpose of a vector of operators. Also, we assume the following commutation relations hold for the noise components:
\begin{equation}
\label{noisecom}
\left[{d\tilde w(t),d\tilde w(t)^{\dagger}} \right] \triangleq d\tilde
w (t)d\tilde w(t)^{\dagger} - (d\tilde w(t)^* d\tilde
w(t)^T)^T = T_w dt.
\end{equation}
Here $T_w$ is a Hermitian matrix and the notation $^T$ denotes the transpose of a vector or matrix of operators. The noise processes can be represented as operators on an appropriate Fock
space; for more details, see \cite{Parthasarathy92}.

The process $\beta_{w}(t)$ represents variables of other systems
which may be passed to the system (\ref{sys_a}) via an interaction. Therefore,
it is required that $\beta_{w}(0)$ be an operator on a Hilbert
space distinct from that of $a_0$ and the noise processes. We also
assume that $\beta_{w}(t)$ commutes with $a(t)$ for all $t\geq
0$. Moreover, since $\beta_{w}(t)$ is an adapted process, we note that $\beta_{w}(t)$ also
commutes with $d\tilde{w}(t)$ for all $t\geq 0$. We have also that $\left[a, dW^\dagger \right]=0$.

\section{Commutation Relations}
For the nonlinear quantum system (\ref{sys_a}), the initial system variables $a(0)=a_{0}$ consist of operators satisfying the commutation relations
\begin{eqnarray}\nonumber
\label{relations}
\left[ {a_j (0),a_k^{*}(0)} \right] &=& \Theta _{jk} ,\qquad j,k = 1,
\ldots ,n.\\
\left[ {\bar a_j (0),\bar a_k^{*}(0)} \right] &=& \bar \Theta _{jk} ,\qquad j,k = 1,
\ldots ,2n.
\end{eqnarray}
Here, the commutator is defined by $\left[{a_j,a_k^{*}} \right]
\triangleq a_ja_k^{*}  - a_{k}^* a_j = \Theta_{jk}$ where $\Theta$ is a complex matrix with elements $\Theta_{jk}$.
With $a^T = (a_1,\ldots, a_n)$, the  relations (\ref{relations}) can  be written as
\begin{eqnarray}\nonumber \label{thetaandbar}
\left[ {a,a^\dagger } \right] & \triangleq & aa^\dagger  - (a^*a^T
)^T= \Theta\\
\left[\bar a, \bar a^\dagger\right]& \triangleq & \left[ {\begin{array}{*{20}{c}}
\Theta &0\\
0&{{\Theta ^*}}
\end{array}} \right]=\bar \Theta.
\end{eqnarray}
\subsection{Preservation of the Commutation Relations}
The following theorem provides an algebraic characterization of
when the nonlinear quantum system (\ref{sys_a}) preserves the commutation
relations (\ref{relations}) as time evolves.
\begin{theorem}\label{thm1}
For the nonlinear quantum system (\ref{sys_a}),
we have that $[\bar a_{p}(0), \bar a_{q}^{*}(0)]=\bar \Theta_{pq}$ implies
$[\bar a_p(t),\bar a_q^{*}(t)]=\bar \Theta_{pq}$ for all $t\geq 0$ with
$p,q=1\ldots, 2n$ if and only if
\begin{eqnarray}\nonumber \label{preservation}
&& \left[ \bar A(\bar a, \bar a^\dagger),\bar a^\dagger \right]+\left[\bar a, \bar A(\bar a,\bar a^\dagger)^\dagger \right]\\\nonumber && \quad + \bar B(\bar a, \bar a^\dagger)T_{ \bar w}\bar B(\bar a,\bar a^\dagger)^\dagger =0;\\\nonumber
&& \left[ \bar B(\bar a,\bar a^\dagger),\bar a^\dagger\right] = 0 \qquad \mbox{and}\\
&& \left[ \bar a, \bar B(\bar a,\bar a^\dagger)^\dagger \right]=0
\end{eqnarray}
with $T_{\bar w} = \left[ {\begin{array}{*{20}{c}}
{T_w}&0\\
0&{T_w^*}
\end{array}} \right]$.
\end{theorem}
\section{Physical Realizability And The Nonlinear Quantum Harmonic Oscillator}
In this section, we consider conditions under which a nonlinear quantum system of the form (\ref{sys_a}) corresponds to a nonlinear open quantum harmonic oscillator. Such a system is said to be physically realizable.
The class of nonlinear open quantum harmonic oscillators under consideration are defined by the nonlinear coupling operator $ \bar L $ and a nonlinear Hamiltonian $ \bar H$.  To derive a nonlinear system of the form (\ref{sys_a}) from a nonlinear open quantum harmonic oscillator defined by $\bar L$ and $\bar H$, we proceed as in Section II to get the following definition.
\begin{definition}\label{def1}
A nonlinear quantum system of the form (\ref{sys_a}) is said to be physically realizable if it is a representation of a nonlinear open quantum harmonic oscillator defined by a nonlinear coupling operator $\bar L$ and a nonlinear Hamiltonian $\bar H$, i.e., if there exist $\bar H$ and $\bar L$ such that the matrices $\bar A(a,a^\dagger), \bar B(a,a^\dagger), \bar C(a,a^\dagger)$ and $\bar D(a,a^\dagger)$ satisfy the following equations
\begin{eqnarray}\label{hamphys}\nonumber
\bar A(a,a^\dagger)&=&\frac{1}{2}\left[\bar L^\dagger, \bar a\right]\bar I \bar L + i\left[\bar H, \bar a\right] ;\\\nonumber
\bar B(a,a^\dagger)&=&\left[\bar L^\dagger, \bar a\right]\bar I;\\\nonumber
\bar C(a,a^\dagger)&=&\bar L;\\
\bar D(a,a^\dagger)&=&I;
\end{eqnarray}
where $\bar I =\left[ {\begin{array}{*{20}{c}}
I&{0}\\
{0}&-I
\end{array}} \right].$
\end{definition}
The following theorem provides necessary and sufficient conditions for a nonlinear quantum system of the form (\ref{sys_a}) to be physically realizable.
\begin{theorem}\label{thm2}
A nonlinear quantum system of the form (\ref{sys_a}) is  physically realizable with $T_{\bar w}= \bar I=\left[ {\begin{array}{*{20}{c}}
{I}&0\\
0&{-I}
\end{array}} \right] $ if and only if
\begin{eqnarray}\nonumber \label{phys_cond1}
&&\left[ \bar A( a,  a^\dagger),\bar a^\dagger \right]+\left[\bar a, \bar A( a, a^\dagger)^\dagger \right]\\\nonumber && \quad + \bar B( a,  a^\dagger)\bar I \bar B( a, a^\dagger)^\dagger =0;\\\nonumber
&& \left[ \bar B( a, a^\dagger),\bar a^\dagger \right] = 0;\\\nonumber
&& \left[ \bar a, \bar B( a, a^\dagger)^\dagger \right]=0;\\\nonumber
&& \bar B( a, a^\dagger)= \left[\bar C(a,a^\dagger)^\dagger,\bar a \right]\bar I\qquad \mbox{and}\\
&&\bar D( a,  a^\dagger)=I.
\end{eqnarray}
In this case, the corresponding nonlinear self-adjoint Hamiltonian $\bar H$ is given by
\begin{equation}\label{Ham}
\bar H=\frac{i}{2\bar n}\left({\bar a^\dag }{\bar \theta ^{ - 1}}\bar A(a,{a^\dag }) - \bar A{(a,{a^\dag })^\dag }{\bar \theta ^{ - 1}}\bar a\right)
\end{equation}
where $\bar n$ is as defined in Definition \ref{def4} and the corresponding coupling operator is $\bar L=\bar C(a,a^\dagger)$.
\end{theorem}
\emph{\textbf{Note:}} Note that the physical realizability conditions as illustrated by conditions (\ref{phys_cond1}) impose the restriction of requiring $C(a)$ to be a vector of first order polynomials of the annihilation operator $a$; i.e, the index $m_{kv}$ in Definition \ref{def4} is only $1$ and $
C_v(a)= \sum\limits_{p = 1}^n \beta_{vp}a_p$
for $v=1,\ldots, m$. Consequently, this in turn imposes the restriction of requiring the matrix $B(a,a^\dagger)$ to be a constant matrix.

\section{Lossless Property of The Class of Nonlinear Quantum Systems}
In this section, we consider a class of nonlinear quantum systems of the form (\ref{sys_a}) and we give a quantum version of the lossless property for this class.
\begin{definition}\label{def2}
Given an operator valued quadratic form
\begin{equation}\label{rate}
\bar r(\bar\beta_w, \bar y)=\bar y^\dagger \bar Q \bar y+ 2\bar y^\dagger \bar S \bar \beta_w+ \bar \beta_w^\dagger \bar R \bar \beta_w
\end{equation}
where $\bar Q \in \mathbb{C}^{\bar p\times \bar p}$, $\bar S \in \mathbb{C}^{\bar p\times \bar m}$ and $\bar R \in \mathbb{\bar D}^{\bar m\times \bar m}$ are constant matrices with $\bar Q$ and $\bar R$ Hermitian, we say that the system (\ref{sys_a}) is dissipative with supply rate $\bar r(\bar \beta_w, \bar y)$ if there exists a storage function $\bar \phi( \bar a)$ which is polynomial operator valued mapping with $\bar \phi(\bar a)\geq0$ for all $t$ such that $\bar \phi(0)=0$ and
\begin{equation}\label{dissipativedef}
\int_0^t {\bar r({\bar \beta _w(s)},\bar y(s))} ds \geq \bar \phi(\bar a)
\end{equation}
for all $t\geq 0$ and all operator valued $\bar \beta_w (t)$.
\end{definition}

\begin{definition}\label{def3}
The nonlinear quantum system of the form (\ref{sys_a}) is said to be lossless if it is dissipative with $\bar Q=-I$, $\bar S=0$ and $\bar R=I$ and
\begin{equation}\label{losslessdef}
\int_0^t {\bar r({\bar \beta _w(s)},\bar y(s))} ds =0
\end{equation}
whenever $\bar \phi(0)=\bar \phi(\bar a(t))=0$ for all operator valued $\bar \beta_w (t)$ and all $t\geq 0$.
\end{definition}
\begin{theorem}\label{thm3}
A necessary and sufficient condition for the nonlinear quantum system (\ref{sys_a}) to be lossless with respect to the supply rate (\ref{rate}) is that there exists a storage function $\bar \phi(\bar a)$ with $\bar \phi(\bar a)\geq0$ for all $t$ such that $\bar \phi(0)=0$ and
\begin{eqnarray}\label{lossless_cond}\nonumber
\nabla \bar \phi {(\bar a)^\dagger }\bar A(a,a^\dagger) &=&  - \bar C{(a,a^\dagger)^\dagger }\bar C(a,a^\dagger)\\\nonumber
\frac{1}{2}\bar B{(a,a^\dagger)^\dagger }\nabla \bar \phi (\bar a) &=&  - \bar C(a,a^\dagger)\\
I-\bar D(a,a^\dagger)^\dagger \bar D(a,a^\dagger)&=&0.
\end{eqnarray}
\end{theorem}
\section{Physical Realizability And The Lossless Property}
In this section, we show that the nonlinear quantum system of the form (\ref{sys_a}) is physically realizable if and only if it is lossless with a specific storage function.
\begin{theorem}\label{thm4}
The nonlinear quantum system (\ref{sys_a}) is physically realizable if and only if it is lossless with a differentiable storage function $\bar \phi(\bar a)$ satisfying
\begin{eqnarray}\label{storage_cond}
\left[\nabla \bar \phi(\bar a),\bar a^T\right]&=& \left[ {\begin{array}{*{20}{c}}
0&{2I}\\
{-2I}&0
\end{array}} \right].
\end{eqnarray}
\end{theorem}
\section{Example}
We consider a nonlinear quantum system described by
\begin{eqnarray}\nonumber
A(a,a^\dagger)&=& \left[ {\begin{array}{*{20}{c}}
{ - {k_{{a_1}}}{a_1} + 2a_1^*a_2^2}\\
{ - {k_{{a_2}}}{a_2} - 2a_2^*a_1^2}
\end{array}} \right];\\\nonumber
B(a,a^\dagger)&=& \left[ {\begin{array}{*{20}{c}}
{ - \sqrt {2{k_{{a_1}}}} }&0\\
0&{ - \sqrt {2{k_{{a_2}}}} }
\end{array}} \right];\\\nonumber
C(a,a^\dagger)&=& \left[ {\begin{array}{*{20}{c}}
{\sqrt {2{k_{{a_1}}}} }a_1\\
{\sqrt {2{k_{{a_2}}}} }a_2
\end{array}} \right];\\\nonumber
D(a,a^\dagger)&=& I.
\end{eqnarray}
where $a$ is the annihilation operator, $a^\dagger$ is the creation operator and $k_{a_1}$ and $k_{a_2}$ are constants.

Hence, \begin{eqnarray}\nonumber
\bar A(a,a^\dagger)&=& \left[ {\begin{array}{*{20}{c}}
{ - {k_{{a_1}}}{a_1} + 2a_1^*a_2^2}\\
{ - {k_{{a_2}}}{a_2} - 2a_2^*a_1^2}\\
{ - {k_{{a_1}}}{a_1^*} + 2a_1a_2^{2*}}\\
{ - {k_{{a_2}}}{a_2^*} - 2a_2a_1^{2*}}
\end{array}} \right];\\\nonumber
\bar B(a,a^\dagger)&=& \left[ {\begin{array}{*{20}{c}}
{ - \sqrt {2{k_{{a_1}}}} }&0&0&0\\
0&{ - \sqrt {2{k_{{a_2}}}} }& 0 & 0\\
0& 0 &{  -\sqrt {2{k_{{a_1}}}} }&0\\
0&0& 0 &{  -\sqrt {2{k_{{a_2}}}} }
\end{array}} \right];\\\nonumber
\bar C(a,a^\dagger)&=& \left[ {\begin{array}{*{20}{c}}
{\sqrt {2{k_{{a_1}}}} }a_1\\
{\sqrt {2{k_{{a_2}}}} }a_2\\
{\sqrt {2{k_{{a_1}}}} }a_1^*\\
{\sqrt {2{k_{{a_2}}}} }a_2^*
\end{array}} \right];\\\nonumber
\bar D(a,a^\dagger)&=& I.
\end{eqnarray}

The matrices $A(a,a^\dagger)$, $B(a,a^\dagger)$ and $C(a,a^\dagger)$ satisfy the properties of the class of nonlinear quantum systems under consideration as given in Definition \ref{def4}. In fact, $\left[B(a,a^\dagger),a^T\right]=0$, $\left[C(a), a^T\right]=0$ and $\left[A(a,a^\dagger), a^T\right]=-\left[a, A(a,a^\dagger)^T\right]= \left[ {\begin{array}{*{20}{c}}
{ - 2a_2^2}&0\\
0&{2a_1^2}
\end{array}} \right]$.

Moreover, $A(a,a^\dagger)$ and $C(a,a^\dagger)$ satisfy conditions (\ref{prop_matrices}) in Definition \ref{def4}. In this case, $\bar n =4$.

Also, condition (\ref{prop_matrices_2}) in Definition \ref{def4} is satisfied. In fact, $\frac{1}{{{2{\bar n}}}}\left[ {\bar A{{(a,{a^\dag })}^\dag }{{\bar \theta }^{ - 1}}\bar a,\bar a} \right]= \left[ {\begin{array}{*{20}{c}}
{a_1^*a_2^2}\\
{ - a_2^*a_1^2}\\
{a_2^{2*}{a_1}}\\
{ - a_1^{2*}{a_2}}
\end{array}} \right]$ and $\frac{1}{{{2{\bar n}}}}\left[ {{{\bar a}^\dag }{{\bar \theta }^{ - 1}}\bar A(a,{a^\dag }),\bar a} \right]=\left[ {\begin{array}{*{20}{c}}
{-a_1^*a_2^2}\\
{  a_2^*a_1^2}\\
{-a_2^{2*}{a_1}}\\
{  a_1^{2*}{a_2}}
\end{array}} \right]$. On the other hand, $\bar A(a,a^\dagger)- \frac{1}{2}\bar B(a,{a^\dag })\bar C(a,a^\dagger)= \left[ {\begin{array}{*{20}{c}}
{2 a_1^*a_2^2}\\
{  -2a_2^*a_1^2}\\
{2 a_2^{2*}{a_1}}\\
{  -2a_1^{2*}{a_2}}
\end{array}} \right]$. Hence, \begin{eqnarray}\nonumber
&&\frac{1}{{{2{\bar n}}}}\left[ {\bar A{{(a,{a^\dag })}^\dag }{{\bar \theta }^{ - 1}}\bar a,\bar a} \right] - \frac{1}{{{2{\bar n}}}}\left[ {{{\bar a}^\dag }{{\bar \theta }^{ - 1}}\bar A(a,{a^\dag }),\bar a} \right]=\\\nonumber&& \qquad \bar A(a,{a^\dag }) - \frac{1}{2}\bar B(a,{a^\dag })\bar C(a,a^\dagger).
\end{eqnarray}
In addition, the system under consideration is physically realizable. In fact, $\bar B( a, a^\dagger)= \left[\bar C(a,a^\dagger)^\dagger,\bar a \right]\bar I= \left[ {\begin{array}{*{20}{c}}
{ - \sqrt {2{k_{{a_1}}}} }&0&0&0\\
0&{ - \sqrt {2{k_{{a_2}}}} }&0&0\\
0&0&{-\sqrt {2{k_{{a_1}}}} }&0\\
0&0&0&{-\sqrt {2{k_{{a_2}}}} }
\end{array}} \right]$.

On the other hand, $\left[ \bar B( a, a^\dagger),\bar a^\dagger \right] = 0$;
$\left[ \bar a, \bar B( a, a^\dagger)^\dagger \right]=0$ and
$\bar D( a,  a^\dagger)=I$.

Also, $\left[ {\bar A(a,{a^\dag }),{{\bar a}^\dag }} \right] = \left[ {\begin{array}{*{20}{c}}
{ - {k_{{a_1}}}}&{4a_1^*{a_2}}&{ - 2a_2^2}&0\\
{ - 4a_2^*{a_1}}&{ - {k_{{a_2}}}}&0&{2a_1^2}\\
{2a_2^{2*}}&0&{{k_{{a_1}}}}&{ - 4{a_1}a_2^*}\\
0&{ - 2a_1^{2*}}&{4{a_2}a_1^*}&{{k_{{a_2}}}}
\end{array}} \right]$,

$\left[ {\bar a,\bar A{{(a,{a^\dag })}^\dag }} \right] = \left[ {\begin{array}{*{20}{c}}
{ - {k_{{a_1}}}}&{ - 4a_1^*{a_2}}&{2a_2^2}&0\\
{4a_2^*{a_1}}&{ - {k_{{a_2}}}}&0&{ - 2a_1^2}\\
{ - 2a_2^{2*}}&0&{{k_{{a_1}}}}&{4{a_1}a_2^*}\\
0&{2a_1^{2*}}&{ - 4{a_2}a_1^*}&{{k_{{a_2}}}}
\end{array}} \right]$ and

$\bar B(a,{a^\dag })\bar I\bar B{(a,{a^\dag })^\dag } = \left[ {\begin{array}{*{20}{c}}
{2{k_{{a_1}}}}&0&0&0\\
0&{2{k_{{a_2}}}}&0&0\\
0&0&{ - 2{k_{{a_1}}}}&0\\
0&0&0&{ - 2{k_{{a_2}}}}
\end{array}} \right]$.

Hence, $\left[ \bar A( a,  a^\dagger),\bar a^\dagger \right]+\left[\bar a, \bar A( a, a^\dagger)^\dagger \right]+ \bar B( a,  a^\dagger)\bar I \bar B( a, a^\dagger)^\dagger =0.$

In this case, the corresponding nonlinear self-adjoint Hamiltonian $\bar H$ is given by
\begin{eqnarray}\nonumber
\bar H&=&\frac{i}{2\bar n}\left({\bar a^\dag }{\bar \theta ^{ - 1}}\bar A(a,{a^\dag }) - \bar A{(a,{a^\dag })^\dag }{\bar \theta ^{ - 1}}\bar a\right)\\\nonumber &=& ia_1^{2*}a_2^2-ia_2^{2*}a_1^2.
\end{eqnarray}
and the corresponding coupling operator is $\bar L=\bar C(a,a^\dagger)$.

Moreover, the system is lossless. In fact, let $\bar \phi(\bar a)=\bar \phi(\bar a)^\dagger=2a_1^*a_1+2a_2^*a_2 \geq 0$. Hence, $\nabla \bar \phi(\bar a)= \left[ {\begin{array}{*{20}{c}}
{2{a_1}}\\
{2{a_2}}\\

{2a_1^*}\\
{2a_2^*}
\end{array}} \right]$. Hence, $\left[\nabla \bar \phi(\bar a),\bar a^T\right]=\left[ {\begin{array}{*{20}{c}}
0&{2I}\\
{-2I}&0
\end{array}} \right]$.

Also, $\nabla \bar \phi {(\bar a)^\dagger }\bar A(a,a^\dagger) =  - \bar C{(a,a^\dagger)^\dagger }\bar C(a,a^\dagger)= -2k_{a_1}(a_1a_1^*+ a_1^*a_1)-2k_{a_2}(a_2a_2^*+ a_2^*a_2)$.

Moreover, $\frac{1}{2}\bar B{(a,a^\dagger)^\dagger }\nabla \bar \phi (\bar a) = - \bar C(a,a^\dagger)$ and $I-\bar D(a,a^\dagger)^\dagger \bar D(a,a^\dagger)=0$.

Hence, the system is lossless using Theorem \ref{thm3}.
\section{Conclusion}
In this paper, the issue of physical realizability for a class of nonlinear quantum systems is developed and was found to be connected to the lossless properties of the system. This provides a connection between the quantum mechanical and nonlinear system theory properties of the system.
%\bibliographystyle{IEEEtran}
%\bibliography{physnonlinearbib}

\begin{thebibliography}{10}
\providecommand{\url}[1]{#1}
\csname url@rmstyle\endcsname
\providecommand{\newblock}{\relax}
\providecommand{\bibinfo}[2]{#2}
\providecommand\BIBentrySTDinterwordspacing{\spaceskip=0pt\relax}
\providecommand\BIBentryALTinterwordstretchfactor{4}
\providecommand\BIBentryALTinterwordspacing{\spaceskip=\fontdimen2\font plus
\BIBentryALTinterwordstretchfactor\fontdimen3\font minus
  \fontdimen4\font\relax}
\providecommand\BIBforeignlanguage[2]{{%
\expandafter\ifx\csname l@#1\endcsname\relax
\typeout{** WARNING: IEEEtran.bst: No hyphenation pattern has been}%
\typeout{** loaded for the language `#1'. Using the pattern for}%
\typeout{** the default language instead.}%
\else
\language=\csname l@#1\endcsname
\fi
#2}}

\bibitem{Maalouf11-2}
A.~I. Maalouf and I.~R. Petersen, ``Bounded real properties for a class of annihilation-operator linear   quantum systems,''
  \emph{IEEE Transactions \mbox{o}n Automatic Control}, vol.~56, no.~4, pp. 786 -- 801, 2011.

\bibitem{James08}
M.~R. James, H.~I. Nurdin, and I.~R. Petersen, ``\mbox{$H^{\infty}$} control of
  linear quantum stochastic systems,'' \emph{IEEE Transactions \mbox{o}n
  Automatic Control}, vol.~53, no.~8, pp. 1787--1803, 2008.

\bibitem{Gardiner85}
C. ~W. Gardiner and M. ~J. Collett, ``Input and output in damped quantum systems: Quantum stochastic differential equations and the master equation,'' \emph{Physical Review A}, vol.~31, no.~6, pp. 3761--3774, 1985.

\bibitem{Bachor04}
H. ~A. Bachor and T.~ C. Ralph, \emph{A Guide to Experiments in Quantum Optics}.\hskip 1em plus 0.5em
  minus 0.4em\relax Weinheim, Germany:Wiley-VCH, 2004.

\bibitem{Parthasarathy92}
K.~R. Parthasarathy, \emph{An Introduction to Quantum Stochastic
  Calculus}.\hskip 1em plus 0.5em minus 0.4em\relax Berlin: Birkhauser, 1992.

\bibitem{Hudson84}
R.~ L. Hudson and K.~ R. Parthasarathy, ``Quantum $\mbox{Ito's}$ formula and stochastic evolutions,'' \emph{Communications in Mathematical Physics}, vol.~93,  pp.
  301--323, 1984.

\bibitem{Walls08}
D. ~F. Walls and G. ~J . Milburn, \emph{Quantum optics}.\hskip
  1em plus 0.5em minus 0.4em\relax Berlin: Springer-Verlag, 2000.

\bibitem{Wiseman10}
H. ~M. Wiseman and G.~ J. Milburn, \emph{Quantum Measurement and Control}.\hskip 1em plus 0.5em minus 0.4em\relax Cambridge University Press:
  UK, 2010.


\bibitem{Gardiner00}
C. ~W. Gardiner and P. Zoller, \emph{Quantum noise}.\hskip 1em plus 0.5em minus 0.4em\relax Berlin:
  Springer, 2000.

\bibitem{James10}
M. ~R. James and J.~ E. Gough, ``Quantum Dissipative Systems and Feedback Control Design by Interconnection,'' \emph{IEEE Transactions \mbox{o}n Automatic Control}, vol.~55, no.~8,  pp.
 1806--1821, 2010.


\bibitem{Shaiju09}
A.~ J. Shaiju and I.~ R. Petersen, ``On the Physical Realizability of General Linear Quantum Stochastic Differential Equations with Complex Coefficients,'' \emph{Proceedings of the joint 48th IEEE Conference on Decision and Control and the 28th Chinese Control Conference, Shanghai, P.R. China, December 16-18, 2009}, pp.
 1422--1427, 2009.

 \bibitem{Maalouf09}
A.~ I. Maalouf and I.~ R. Petersen, ``Coherent $\mbox{H}^\infty$ Control for a Class of Linear Complex Quantum Systems,'' \emph{Proceedings of the American Control Conference ACC 2009, St. Louis, Missouri}, 2009.

\bibitem{Maalouf11-1}
A.~I. Maalouf and I.~R. Petersen, ``Coherent \mbox{$H^{\infty}$} Control for a class of annihilation-operator linear   quantum systems,''
  \emph{IEEE Transactions \mbox{o}n Automatic Control}, vol.~56, no.~2, pp. 309-319, 2011.


\bibitem{Hassen10}
S~.Z~.S~. Hassen and I.~ R. Petersen and E.~ H. Huntington and M. Heurs and M. ~R. James, ``$\mbox{LQG}$ control of an optical squeezer,'' \emph{Proceedings of the American Control Conference ACC 2010, Baltimore, USA},   pp.2730 - 2735, 2010.

% \bibitem{Maalouf12}
%A.~ I. Maalouf and I.~ R. Petersen, ``On the Physical Realizability of a Class of Nonlinear Quantum Systems,'' \emph{arXiv:1207.5299 [math.OC]}, 2012.
\end{thebibliography}

\end{document}